\documentclass{gtart_h}

\def\ifplaintex{\expandafter\ifx\csname documentclass\endcsname\relax}

\def\gtp{{\mathsurround=0pt\it $\cal G\mskip-2mu$eometry \&\ 
$\cal T\!\!$opology $\cal P\!$ublications}}  

\def\recd{{\small Received:\qua\receiveddate\ifx\reviseddate\relax
\else\qquad Revised:\qua\reviseddate\fi\par}} 

\def\lognumber#1{\def\thelognumber{#1}}
\def\volumenumber#1{\def\thevolumenumber{#1}}
\def\volumeyear#1{\def\thevolumeyear{#1}}
\def\papernumber#1{\def\thepapernumber{#1}}
\def\pagenumbers#1#2{\def\startpage{#1}\def\finishpage{#2}}
\def\published#1{\def\publishdate{#1}}

\def\received#1{\def\receiveddate{#1}}
\def\revised#1{\def\reviseddate{#1}}
\def\accepted#1{\def\accepteddate{#1}}

\def\asciiaddress#1{\def\theasciiaddress{#1}}
\def\asciiemail#1{\def\theasciiemail{#1}}
\def\asciiurl#1{\def\theasciiurl{#1}}

\long\def\asciiabstract#1{\long\def\theasciiabstract{#1}}
\def\asciikeywords#1{\def\theasciikeywords{#1}}

\let\\\par\let\thelognumber\relax\let\thevolumenumber\relax
\let\thepapernumber\relax\let\thevolumeyear\relax\let\startpage\relax
\let\finishpage\relax\let\publishdate\relax\let\receiveddate\relax
\let\reviseddate\relax\let\accepteddate\relax\let\theasciititle\relax
\let\theasciiauthors\relax\let\theasciiaddress\relax
\let\theasciiabstract\relax\let\theasciikeywords\relax

\let\theasciiemail\relax
\let\theasciiurl\relax

\ifplaintex
\font\logobig=cmssbx10 scaled 3836
\font\logomed=cmssbx10 scaled 2557
\else
\font\logobig=cmssbx10 scaled 4200
\font\logomed=cmssbx10 scaled 2800
\fi

\long\def\makeagttitle{   
\count0=\startpage
\agt\hfill      
\hbox to 45truept{\vbox to 0pt{\vglue -13truept{\logomed A\kern -.37em{\logobig 
T}\kern -.38em G}\vss}\hss}
\break
{\small Volume \thevolumenumber\ (\thevolumeyear)
\startpage--\finishpage\nl
Published: \publishdate}

\vglue .25truein

{\parskip=0pt\leftskip 0pt plus
1fil\def\\{\par\smallskip}{\Large\bf\thetitle}\par\medskip} \vglue
0.05truein

{\parskip=0pt\leftskip 0pt plus 1fil\def\\{\par}{\sc\theauthors}
\par\medskip}%
 
\vglue 0.03truein 

{\small\leftskip 25truept\rightskip 25truept{\bf Abstract}\stdspace\theabstract

{\bf AMS Classification}\stdspace\theprimaryclass
\ifx\thesecondaryclass\relax\else; \thesecondaryclass\fi\par
{\bf Keywords}\stdspace \thekeywords\par}\vglue 7truept

}   

\ifplaintex
\font\phead=cmsl9 scaled 950
\font\pnum=cmbx10 scaled 913
\font\pfoot=cmsl9 scaled 950
\headline{\vbox to 0pt{\vskip -4.5mm\line{\small\phead\ifnum
\count0=\startpage ISSN 1472-2739 (on-line) 1472-2747 (printed)
\hfill {\pnum\folio}\else\ifodd\count0\def\\{ }%
\ifx\theshorttitle\relax\thetitle\else\theshorttitle\fi\hfill{\pnum\folio}
\else\def\\{ and }{\pnum\folio}\hfill\ifx\theshortauthors\relax\theauthors
\else\theshortauthors\fi\fi\fi}\vss}}
\footline{\vbox to 0pt{\vglue 0mm\line{\small\pfoot\ifnum\count0=\startpage
\copyright\ \gtp\hfill\else
\agt, Volume \thevolumenumber\ (\thevolumeyear)\hfill\fi}\vss}}
\else
\font=cmsl9 scaled 1050
\font\lnum=cmbx10 
\font=cmsl9 scaled 1050
\makeatletter
\def\@oddhead{{\small\lhead\ifnum\count0=\startpage ISSN 1472-2739 
(on-line) 1472-2747 (printed)\hfill {\lnum\number\count0}\else\ifodd\count0
\def\\{ }\ifx\theshorttitle\relax \thetitle \else\theshorttitle\fi\hfill
{\lnum\number\count0}\else\def\\{ and }{\lnum\number\count0}
\hfill\ifx\theshortauthors\relax 
\theauthors\else\theshortauthors\fi\fi\fi}}\def\@evenhead{\@oddhead}
\def\@oddfoot{\small\lfoot\ifnum\count0=\startpage\copyright\ \gtp\hfill\else
\agt, Volume \thevolumenumber\ (\thevolumeyear)\hfill\fi}
\def\@evenfoot{\@oddfoot}
\makeatother
\fi
\let\maketitlepage\makeagttitle
\let\makeshorttitle\maketitlepage
\let\maketitle\maketitlepage

\newwrite\gtoutfile
\long\gdef\makeheadfile{  
{\def\\{, }\def\s{ }
\immediate\openout\gtoutfile head.xxx
\immediate\write\gtoutfile{Proxy-for: \ifx\theasciiauthors\relax
\theauthors\else\theasciiauthors\fi\s<\ifx\theasciiemail\relax\theemail\else\theasciiemail\fi>}
\immediate\write\gtoutfile{\noexpand\\}
\immediate\write\gtoutfile{Authors: \ifx\theasciiauthors\relax
\theauthors\else\theasciiauthors\fi}
{\def\\{ }\immediate\write\gtoutfile{Title: \ifx\theasciititle\relax
\thetitle\else\theasciititle\fi}}
\immediate\write\gtoutfile{Subj-class: GT or SG, GR etc}
\immediate\write\gtoutfile{MSC-class: \theprimaryclass\ifx\thesecondaryclass\relax\else, \thesecondaryclass\fi}
\immediate\write\gtoutfile{Journal-ref: Algebr. Geom. Topol. \thevolumenumber\s
(\thevolumeyear) \startpage-\finishpage}
\immediate\write\gtoutfile{Comments: Published by Algebraic and
Geometric Topology at}
\immediate\write\gtoutfile{\s\s\s  http://www.maths.warwick.ac.uk/agt/AGTVol\thevolumenumber/agt-\thevolumenumber-\thepapernumber.abs.html}
\immediate\write\gtoutfile{\noexpand\\}
\immediate\write\gtoutfile{}
\ifx\theasciiabstract\relax
\immediate\write\gtoutfile{\theabstract}\else
\immediate\write\gtoutfile{\theasciiabstract}\fi
\immediate\write\gtoutfile{}
\immediate\write\gtoutfile{\noexpand\\}
\immediate\write\gtoutfile{}
\immediate\closeout\gtoutfile}}  

\def\maketitlepage{\makeagttitle\makeheadfile}
\let\makeshorttitle\maketitlepage
\let\maketitle\maketitlepage

\lognumber{17}
\volumenumber{5}
\volumeyear{2005}
\papernumber{17}
\pagenumbers{379}{403}
\received{16 December 2004} 
\revised{21 April 2005}
\accepted{6 May 2005}
\published{22 May 2005}

\usepackage{amssymb,epsfig,amsmath}
\allowdisplaybreaks

\theoremstyle{plain}

\newtheorem{theorem}{Theorem}
\newtheorem{proposition}{Proposition}[section]

\newtheorem{conjecture}{Conjecture}

\theoremstyle{definition}
\newtheorem{definition}[proposition]{Definition}

\theoremstyle{remark}

\newtheorem{remark}[proposition]{Remark}

\newcommand{\psdraw}[2]
         {\begin{array}{c} \hspace{-1.3mm}
        \raisebox{-4pt}{\epsfig{figure=draws/#1.eps,width=#2}}
        \hspace{-1.9mm}\end{array}}
\newcommand{\ppsdraw}[2]
         {\begin{array}{c} \hspace{-1.3mm}
        \raisebox{-4pt}{\epsfig{figure=draws/#1.eps,height=#2}}
        \hspace{-1.9mm}\end{array}}

\def\BN{\mathbb N}
\def\BZ{\mathbb Z}

\def\BQ{\mathbb Q}
\def\BR{\mathbb R}
\def\BC{\mathbb C}

\def\l{\lambda}

\def\S{\Sigma}

\def\s{\sigma}

\def\la{\langle}
\def\ra{\rangle}

\def\e{\epsilon}

\def\d{\delta}
\def\b{\beta}

\def\s{\sigma}

\def\SL{\mathrm{SL}}

\def\z{\zeta}
\def\VC{\text{VC}}
\def\vol{\mathrm{vol}}

\def\qb#1#2{\left[\begin{matrix} #1 \\ #2 \end{matrix}\right]}

\begin{document}


\title[Experimental evidence for the Volume Conjecture]{Experimental evidence for
the Volume Conjecture\\for
the simplest hyperbolic non-2-bridge knot}

\author{Stavros Garoufalidis\\Yueheng Lan}

\address{School of Mathematics,  Georgia Institute of 
Technology\\Atlanta, GA 30332-0160, 
USA\\{\rm and}\\School of Physics, Georgia Institute of Technology\\Atlanta, 
GA 30332-0160, USA}

\asciiaddress{School of Mathematics,  Georgia Institute of 
Technology\\Atlanta, GA 30332-0160, 
USA\\and\\School of Physics, Georgia Institute of Technology\\Atlanta, 
GA 30332-0160, USA}

\asciiemail{stavros@math.gatech.edu, gte158y@mail.gatech.edu}
\gtemail{\mailto{stavros@math.gatech.edu}, \mailto{gte158y@mail.gatech.edu}}

\gturl{\url{http://www.math.gatech.edu/~stavros}, 
\url{http://cns.physics.gatech.edu/~y-lan}}

\asciiurl{http://www.math.gatech.edu/ stavros, 
http://cns.physics.gatech.edu/ y-lan}

\begin{abstract}
Loosely speaking, the Volume Conjecture states that the limit of the
$n$-th colored Jones polynomial of a hyperbolic knot, evaluated at the
primitive complex $n$-th root of unity is a sequence of complex
numbers that grows exponentially. Moreover, the exponential growth
rate is proportional to the hyperbolic volume of the knot.
 
We provide an efficient formula for the colored Jones function
of the simplest hyperbolic non-2-bridge knot, and using this formula,
we provide numerical evidence for the Hyperbolic Volume
Conjecture for the simplest hyperbolic non-2-bridge knot.
\end{abstract}

\asciiabstract{%
Loosely speaking, the Volume Conjecture states that the limit of the
n-th colored Jones polynomial of a hyperbolic knot, evaluated at the
primitive complex n-th root of unity is a sequence of complex numbers
that grows exponentially. Moreover, the exponential growth rate is
proportional to the hyperbolic volume of the knot.  We provide an
efficient formula for the colored Jones function of the simplest
hyperbolic non-2-bridge knot, and using this formula, we provide
numerical evidence for the Hyperbolic Volume Conjecture for the
simplest hyperbolic non-2-bridge knot.}

\asciikeywords{Knots, q-difference equations, asymptotics, 
Jones polynomial, Hyperbolic Volume Conjecture, character varieties,
recursion relations, Kauffman bracket, skein module, fusion, SnapPea, m082}

\maketitle


\section{Introduction}
\label{sec.intro}

\subsection{The Hyperbolic Volume Conjecture}
\label{sub.HVC}

The {\em Volume Conjecture} connects two very
different approaches to knot theory,  
namely Topological Quantum Field Theory and Riemannian
(mostly Hyperbolic) Geometry. The Volume Conjecture states that 
for every hyperbolic knot $K$ in $S^3$ 
\begin{equation}
\label{eq:vc}
\lim_{n \to \infty} \frac{\log|J_K(n)(e^{\frac{2 \pi i}{n}})|}{n}=
\frac{1}{2 \pi}\, \vol(K),
\end{equation}
where 
\begin{itemize}
\item
$J_K(n) \in \BZ[q^{\pm}]$ is the Jones polynomial of a knot colored with the 
$n$-dimensional irreducible representation of $\mathfrak{sl}_2$,
normalized so that it equals to $1$ for the unknot (see \cite{J, Tu}), and 
\item 
$\vol(K)$ is the {\em volume} of a complete hyperbolic metric 
in the knot complement $S^3-K$; see \cite{Th}. 
\end{itemize}
The conjecture was formulated in this form by Murakami-Murakami \cite{MM} 
following an earlier version due to Kashaev, \cite{K}. 
The Volume Conjecture is an analytic question that contains two parts: 
\begin{itemize}
\item[(a)] to show that
a limit of a sequence of complex numbers exists, and 
\item[(b)] to identify the 
limit with a known geometric invariant of a knot.
\end{itemize}

Currently, the Volume Conjecture 
is known for only for the simplest hyperbolic knot: the
$4_1$ (due to Ekholm), 
\cite{M}. For the $4_1$ knot, due to special circumstances there is
a 1-dimensional sum formula for $J_{4_1}(n)(e^{2 \pi i/n})$ where
the summand is a closed-form {\em positive} term. 
In that case, it is elementary to see that  the volume conjecture holds.  
There is no other known knot that exhibits similar behavior.

Despite the optimism and the belief in the conjecture within the area of
Quantum Topology, there is natural suspicion about it, due to lack of
evidence.

There are several difficulties in the Volume Conjecture for a knot $K$: 
the left hand side involves a sequence of polynomials $J_K(n)$ (for $n=1,2,
\dots$) with little understood relation to the geometry of the knot
complement. In fact, a major problem is to give {\em efficient multisum
formulas} for the polynomials $J_K(n)$. A naive approach based on 
cables computes $J_K(n)$ in $2^{c n^2}$ steps, where $n$ is the number of 
crossings. An alternative state-sum formula (see \cite[Section 3.2]{GL})
computes $J_K(n)$ in $n^c$ steps. In either case, when $c=18$ and $n=500$,
these formulas are inefficient in numerically computing the left hand side
of \eqref{eq:vc}.

For a fixed knot, the polynomials $J_K(n)$ are not random. In \cite{GL},
it was proven that they satisfy a linear $q$-difference equation. Moreover,
it was explained that 
the results of \cite{GL} together with the WZ algorithms (see \cite{WZ})
can in principle compute the above mentioned linear $q$-difference equation.
However computing the $q$-difference equation for knots with $6$ crossings
is already a difficult task. Thus, this method is of little use in 
numerical computations of the left hand side of \eqref{eq:vc}.

In \cite{Ga} and \cite{GG},  the first author announced a program to prove
the existence of the limit in the Volume Conjecture 
using asymptotics of solutions of linear $q$-difference equations. 
The main idea is that
asymptotics of solutions of difference equations with a small parameter are
governed by the average (on the unit circle) of the corresponding eigenvalues.
When the eigenvalues do not collide or vanish, then this analysis was carried
out in \cite{GG}. In the Volume Conjecture, 
it is known that the eigenvalues collide (ven for the $4_1$ knot),
and the analysis will be extended to cover this case in \cite{CG}.
Combined with the AJ
Conjecture of \cite{Ga}, this program might identify the limit with the
hyperbolic volume. 

Among the hyperbolic knots, the 2-bridge knots are a tractable family 
with well-understood representation theory. 

The purpose of the present paper is to provide: 

\begin{itemize}
\item
an efficient formula for the colored Jones function of the 
simplest hyperbolic non 2-bridge knot: the $k4_3$ knot;
\item
using that formula, {\em numerical evidence} for the Volume Conjecture 
for the $k4_3$ knot.
\end{itemize}

The key idea for efficient formulas for $J_K(n)$ is {\em fusion}
as we explain below.

Aside from the above goals, we provide computation of numerous
topological and geometric invariants of $k4_3$, such as a presentation 
of its fundamental group, peripheral system, Alexander polynomial,
$A$-polynomial, rank
of Heegaard Floer Homology, volume of special Dehn fillings, 
invariant trace fields, as well as a multisum formula for the colored 
Jones function.

Let us end this introduction with an observation.
It appears that for the knot $K_0=k4_3$, the sequence that appears on 
the left hand side of Equation \eqref{eq:vc} is 
eventually {\em monotone decreasing}. 
The same behavior is also exhibited for several $2$-bridge knots that
the authors tried. Monotonicity is an important clue for proving that the 
limit in the HVC exists. Let us formulate this as a conjecture:

\begin{conjecture}
\label{conj.decreasing}
For every knot $K$ in $S^3$, the sequence
$$
\frac{\log|J_K(n)(e^{\frac{2 \pi i}{n}})|}{n}
$$
is eventually decreasing, and bounded above by zero.
\end{conjecture}

\subsection{The knot $K_0$}
\label{sub.candidate}

The word ``simplest'' does not refer to the number of crossings, but rather to
the small number of ideal tetrahedra that are needed for an {\em ideal 
triangulation} of its complement.

Consider the
{\em twisted torus knot} $K_0$, obtained from the torus
knot $T(3,8)$ by adding a full positive twist to two 
strands. Actually, $K_0$ is the second simplest hyperbolic 
non-2-bridge knot with $4$ ideal tetrahedra. The simplest hyperbolic 
non-2-bridge knot $k3_1$ is the Pretzel knot $(-2,3,7)$ requiring
$3$ ideal tetrahedra, with same 
volume as the $5_2$ knot; see \cite{CDW}. However, the colored Jones
function of $k3_1$ is rather complicated.

$K_0$ has
braid presentation $b^2 (ab)^8$, where $a=\s_1, b=\s_2$ are the standard
generators of the braid group with $3$ strands.

\begin{figure}[htpb]
$$ \psdraw{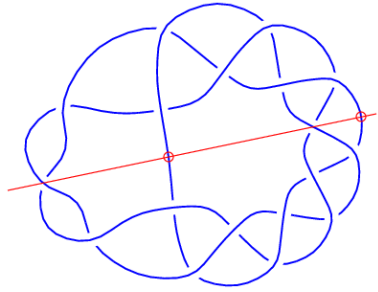}{2in} $$
\caption{The knot $K_0$, and an involution which negates its meridian and 
longitude.}\label{f.k43}
\end{figure}

$K_0$ is a positive hyperbolic knot with $18$ crossings. 
SnapPea identifies $K_0$ with the knot
$k4_3$ of the new census, also known as $m082$ in the old
census.

The notation $k4_3$ reveals that $K_0$ is a hyperbolic knot whose complement 
can be triangulated with $4$ hyperbolic ideal tetrahedra. It is the 
{\em simplest hyperbolic} non-$2$-bridge knot, as was discovered by
\cite{CDW}.

$K_0$ belongs to a family of twisted torus knots, and this family populates the
census of simplest hyperbolic manifolds other than the family of $2$-bridge
knots.

The {\em Dowker code} of $K_0$ is:
$$
18 \ \ 1  \ \  14 \ \  -16 \ \  18 \ \  -20 \ \  22 \ \  -24 \ \  -26 \ \  28 \ \  -30 \ \  32 \ \  -34
\ \  36 \ \  -12 \ \  -2 \ \  4 \ \  -6 \ \  8 \ \  -10
$$
$K_0$ has symmetry group $\BZ/2$. An involution that negates
the meridian and longitude is shown in Figure \ref{f.k43}.

\subsection{The quantum topology of $K_0$}
\label{sub.qt}

Let us define the {\em quantum integer}, the {\em quantum factorial} 
and the {\em quantum binomial coefficients} by:
$$
[n]=\frac{q^{n/2}-q^{-n/2}}{q^{1/2}-q^{-1/2}}
\qquad
[n]!=[1][2]\dots [n] \qquad \qb n k =\frac{[n]!}{[k]![n-k]!},
$$
where $[0]!=1$. The quantum binomial coefficients satisfy
the following recursion relation:
$$
\qb n k=q^{-k/2} \qb {n-1} k + q^{(n-k)/2} \qb {n-1}{k-1},
$$
from which follows that $\qb n k \in \BZ[q^{\pm 1/2}]$.

For natural numbers $a,b$ with $b \leq a$, we denote
$$
[a,b]!=\frac{[a]!}{[b]!},
$$
which also lies in the ring $\BZ[q^{\pm 1/2}]$. 

\begin{definition}
\label{def.Jn}
Let $J_K(n)$ denote the colored Jones polynomial of a knot $K$ in
$S^3$ using the 
$n$-dimensional irreducible representation of $\mathfrak{sl}_2$, 
normalized to equal to $1$ for 
the unknot; see \cite{J,Tu}. 
\end{definition}

The following theorem gives a triple sum formula for $J_{K_0}(n)$,
where the summand is {\em proper $q$-hypergeometric}.

\begin{theorem}
\label{thm.CJ}
For every natural number $n$, we have:
\begin{align*}
&J_{K_0}(n+1)\\ &\qua= \frac{1}{[n+1]} \sum_{k=0,2}^{2n}\sum_{l=|n-k|,2}^{n+k}\sum_z
(-1)^{\frac{k}{2}+z} 
q^{-\frac{3}{8}(2k+k^2)+\frac{7}{8}(2l+l^2)-\frac{51}{8}(2n+n^2)}
\frac{[k+1][l+1]}{[\frac{2n+k}{2}+1]!}  \\
&  \qquad\qquad
\qb {\frac{k+l-n}{2}}{\frac{n+2k+l}{2}-z}
\qb {\frac{n+l-k}{2}}{\frac{3n+l}{2}-z}
\qb {\frac{n+k-l}{2}}{\frac{2n+2k}{2}-z} \\
&  \qquad\qquad 
 \left[\frac{k}{2}\right]!^2 \left[\frac{2n-k}{2},z-\frac{n+k+l}{2}\right]! 
\left[z+1,\frac{n+k+l}{2}+1\right]! \end{align*}
\end{theorem}

Here $\sum_{k=a,2}^b$ means summation for $k=a,a+2,\dots,b$, for even $b-a$.
Although the $z$ summation is infinite, only finitely many terms
contribute. Explicitly, only the terms with $z \in \BN$, 
$$
\mathrm{max}\left(\frac{2n+k}{2},\frac{n+k+l}{2}\right)
\leq z \leq 
\mathrm{min}\left(\frac{n+2k+l}{2},\frac{3n+l}{2},\frac{2n+2k}{2}\right)
$$ 
contribute.

\begin{remark}
\label{rem.qbinom}
The denominators in the above formula come from $[(2n+k)/2+1]!$
and $[n+1]$.
When we evaluate the summand at $q=e^{2 \pi i/(n+1)}$, the order of
the pole is $2$: $1$ from  $[(2n+k)/2+1]!$ and $1$ from $[n+1]$. 
\end{remark}

\begin{remark}
\label{rem.GL}
In \cite[Section 3.2]{GL} using properties of $R$-matrices, Le and the first 
author gave a canonical multisum formula for 
the colored Jones function of a knot
which is presented as the closure of a braid. Specifically, if
the braid has $c$ crossings, then the multisum formula involves summation
over $c$ variables. In the case of the knot $K_0$, it involves summation
over $18$ variables, which makes it intractable, symbolically or numerically.
\end{remark}

\begin{remark}
\label{rem.naive}
The well-known 3-term {\em skein relation} for the Jones polynomial
$$
q \left(\psdraw{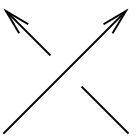}{.3in}\right)-q^{-1} \left(\psdraw{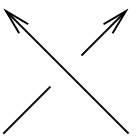}{.3in}\right)
=(q^{1/2}-q^{-1/2}) \left(\psdraw{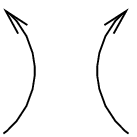}{.3in}\right)
$$
allows one 
to compute naively the Jones polynomial of a knot with $c$ crossings 
in $2^c$ steps. By taking parallels, it allows one to compute the
$n$-th colored Jones polynomial of a knot in $2^{c n^2}$ steps. For
$n=500$, and $c=18$, this requires  
$$
2^{18 \times 500^2} \approx 10^{1354634}
$$
terms. 
\end{remark}

\section{Numerical confirmation of the HVC for $K_0$}
\label{sec.numerical}

Let 
\begin{equation}
\label{eq:HVC}
\VC(n)=\frac{2\pi}{n} \log J_{K_0}(n)|_{q=e^{\frac{2 \pi i}{n}}}.
\end{equation}

Using a {\tt fortran} program by the second author (presented in the 
Appendix), we have computed $\VC(n)$ for $n$ up to $502$.

Here is a table of $\text{Re}(\VC(n))$ versus $n$:
$$
\psdraw{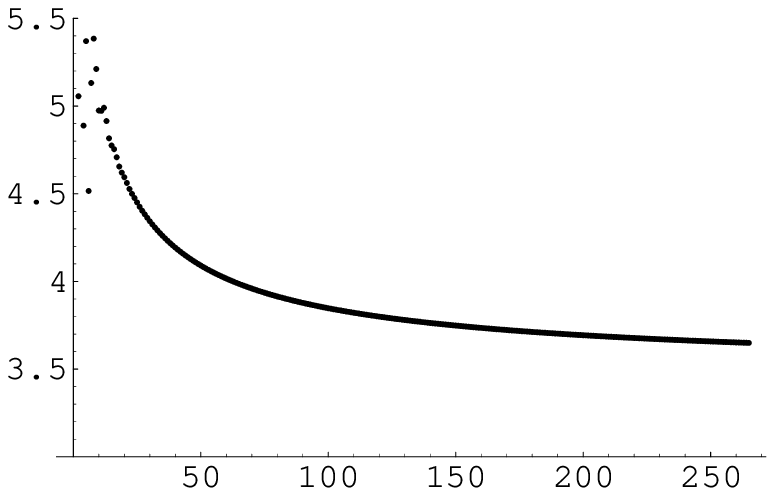}{3in}
$$
Here is a table of $\text{Im}(\VC(n))$ versus $n$:
$$
\psdraw{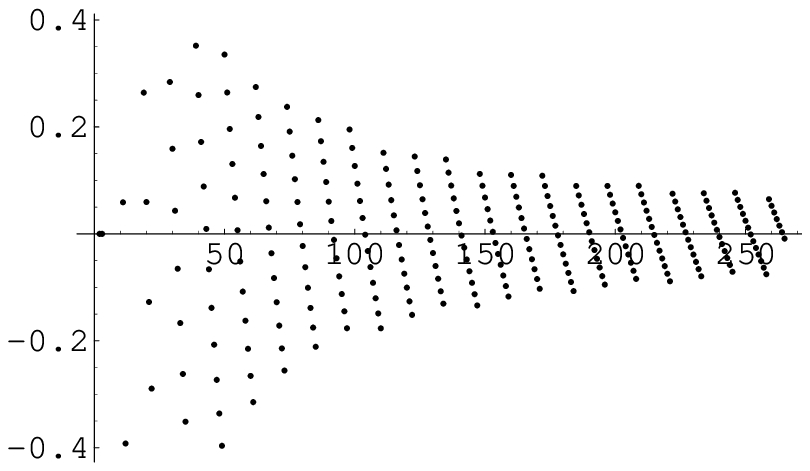}{3in}
$$
Here is a table of $\text{Re}(\VC(n))$ versus $\text{Im}(\VC(n))$:
$$
\psdraw{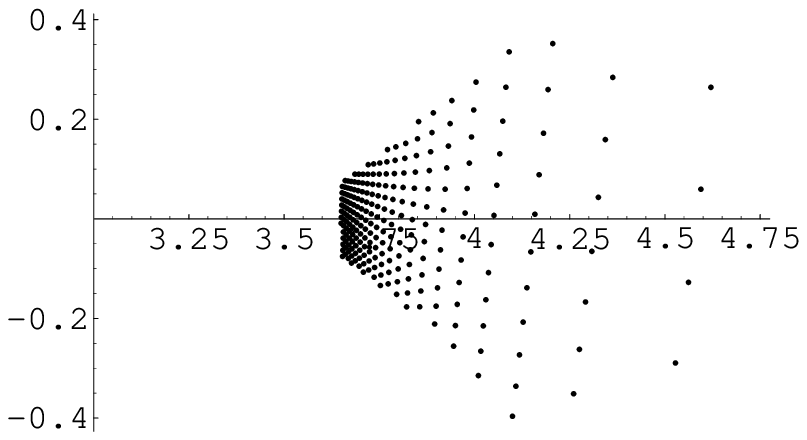}{3in}
$$
There is an evident $12$-fold periodicity that reflects the fact 
that $K_0$ is a twisted $(3,8)$ torus knot.

Quantum Field Theory predicts an asymptotic expansion of the form:
\begin{equation}
\label{eq:assexpansion}
J_{K_0}(n)(e^{2 \pi i/n}) \sim_{n \to \infty}
e^{n \vol(K_0)/(2 \pi) + n i \text{CS}(K_0)} n^{3/2} \left(
C_0 + \frac{C_1}{n} + \frac{C_2}{n^2} + \dots \right)
\end{equation}
for constants $C_i$ (which depend on the knot $K_0$) with $C_0 \neq 0$, 
and for $\text{CS}(K_0) \in \BR/\BZ$. 

Fitting the data with least-squares 
for $n=21,\dots,250$ with $\log(n)/n, \, 1$ and $1/n$ in Mathematica
gives:
{\small
\begin{verbatim}
In[1]:=Fit[data, {Log[n]/n, 1, 1/n}, n]
\end{verbatim}
}
gives 
{\small
\begin{verbatim}
Out[1]:=
3.4750687755045777 - 
    5.518475184459029/n + 9.282495203373793 Log[n]/n
\end{verbatim}
}
and plotting the data against the fit gives:
$$
\psdraw{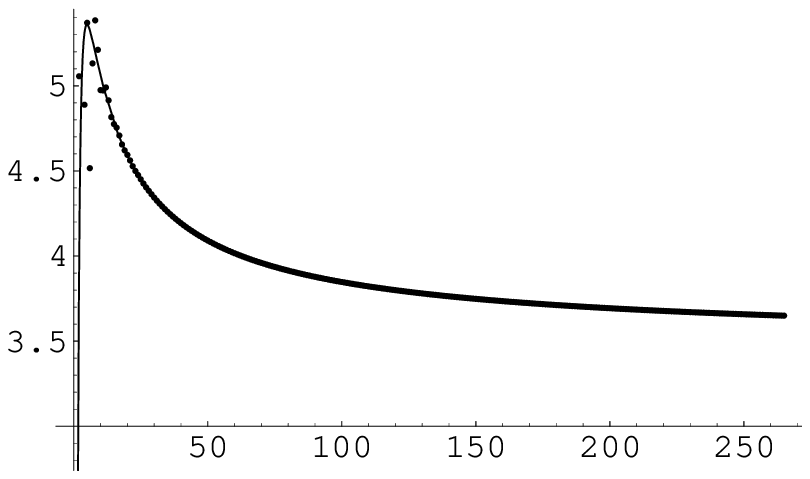}{3in}
$$
Here is a plot of the square of the error:
$$
\psdraw{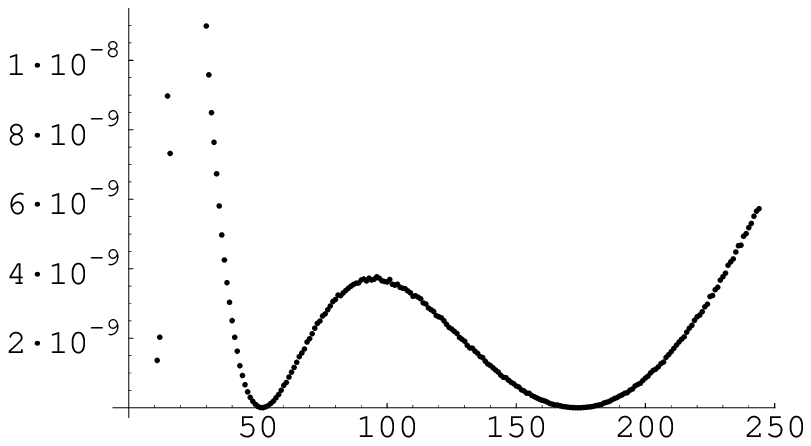}{3in}
$$
This shows that the error of the fit is within $10^{-4}$. 

Notice that $\vol(K_0)=3.474247 \ldots$, 
and $9.282495203373793/(2 \pi)=1.47735 \ldots$. 
Thus, the above data provides strong numerical evidence for the terms
$\vol(K_0)$ and $3/2$ in Equation \eqref{eq:assexpansion}.

Here is a table of $\text{Re}(\VC(2+10n))$ versus $n$ for $n=0,\dots,50$:
$$
\psdraw{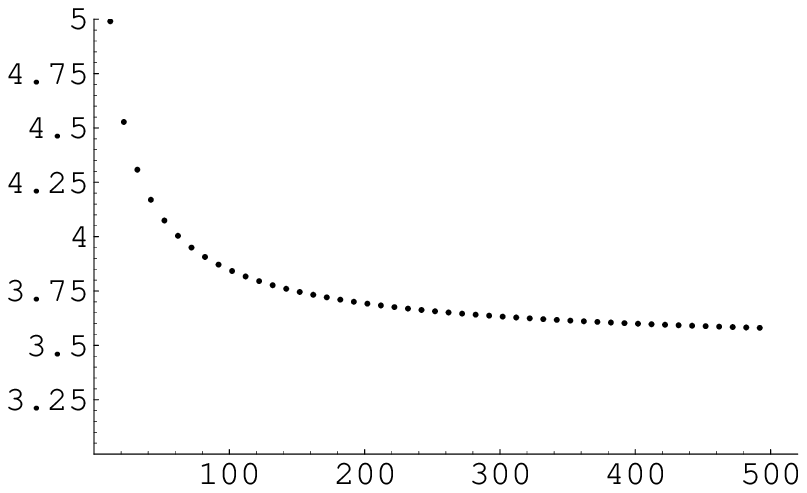}{3in}
$$
Here is a table of $\text{Im}(\VC(2+10n))$ versus $n$ for $n=0,\dots,50$:
$$
\psdraw{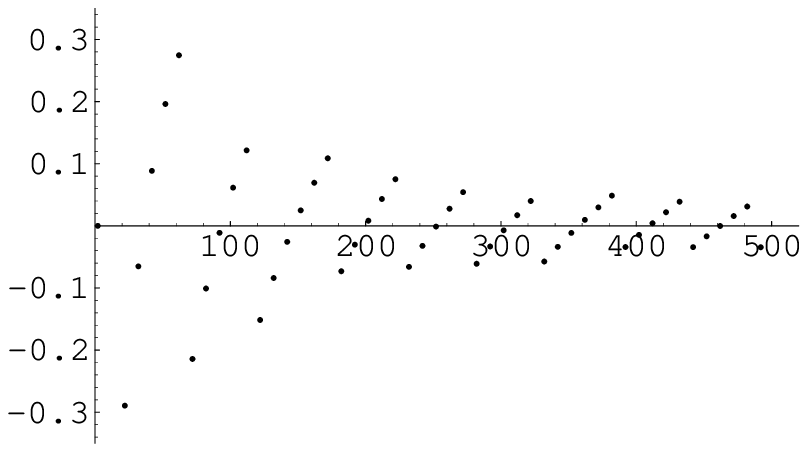}{3in}
$$
Here is a table of $\text{Re}(\VC(n))$ versus $\text{Im}(\VC(2+10n))$
for $n=0,\dots,50$:
$$
\psdraw{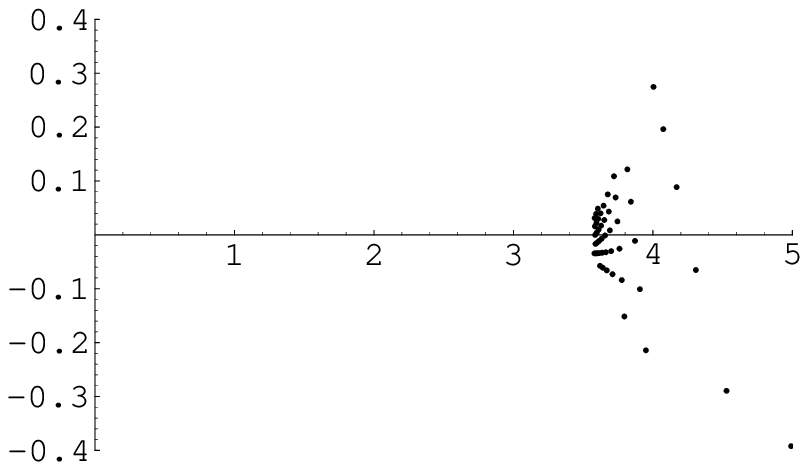}{3in}
$$

\section{Proofs}
\label{sec.proofs}

\subsection{Fusion}
\label{sub.fusion}

The colored Jones function $J_K(n) \in \BZ[q^{\pm}]$ of a knot $K$ in $S^3$
was originally defined using $R$-{\em matrices}, that is solutions to the 
{\em Yang-Baxter equations}; see \cite{J} and \cite{Tu}. These solutions 
are intimately related to the representation theory of quantum groups.
In \cite{GL} Le and the first author used the theory of $R$-matrices to
give formulas for the colored Jones function of a knot $K$. Given a braid
presentation $\b$ of a knot $K$ (that is, a word of length $c$ in the braid
group of some number of strands), one constructs a proper $q$-hypergeometric
function $F_{\b}(n,k_1,\dots,k_c)$ such that
$$
J_K(n)=\sum_{k_1,\dots,k_c=0}^\infty F_{\b}(n,k_1,\dots,k_c).
$$ 
The good thing is that for fixed $n \in \BN$, only finitely many terms
are nonzero. Moreover, the function $F_{\b}$ takes values in the ring
$\BZ[q^{\pm 1/2}]$, i.e., it has no denominators. Thus, we can commute
summation and evaluation at an complex $n$-th root of unity.

The bad thing is that the above sum is multidimensional, which makes
evaluation, symbolic or numerical, impractical.

Thus, we have to find alternative form of presenting $J_K(n)$. To achieve
this, we will use the Kauffman bracket skein module, and  
{\em fusion}. The latter molds together pieces of a knot, producing 
knotted trivalent graphs. An advantage of fusion is that it deals nicely
with twists (see the figures below), and thus it cuts down on the number
of summation variables. A disadvantage of fusion is that it produces 
denominators which vanish when we evaluate at complex $n$th roots of unity.
So, fusion gains low number of summation variables, at the cost of
producing denominators. Of course, if one expands out all the terms, then
the denominators will eventually cancel out, since the colored Jones function
is a Laurent polynomial. However, expanding out denominators and canceling
is too costly, and not efficient enough. The programs developed by the 
second author demonstrate that we can control fusion numerically, 
even when evaluating at complex roots of unity.

\subsection{Review of the Kauffman bracket of knotted trivalent graphs}
\label{sub.review}

Let us recall what is fusion.
A standard reference for this section is \cite{MV} and \cite{KL}. 
We will work with the Kauffman bracket
skein module of links, summarized in the following figure:
$$ 
\psdraw{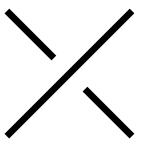}{0.3in}=A \psdraw{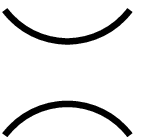}{0.3in}+A^{-1} \psdraw{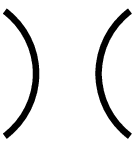}{0.3in},
\hspace{1cm}
\psdraw{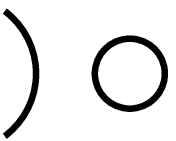}{0.4in} =-(A^2+A^{-2}) \psdraw{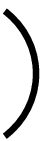}{0.1in}
$$

Here and below, $A=q^{1/4}$. 
The proof of Theorem \ref{thm.CJ}
will use {\em fusion}, leading to {\em knotted trivalent graphs};
the latter are framed embedded colored trivalent graphs. Let us recall
the basic rules of fusion, from \cite{MV} (see also \cite{KL}):
$$
\ppsdraw{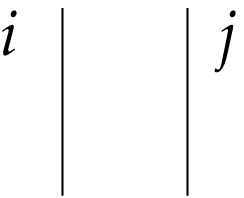}{1.3cm}=\sum_k \frac{\la k\ra}{\la i,j,k \ra} 
\ppsdraw{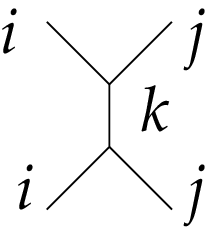}{1.4cm}
$$
$$
\ppsdraw{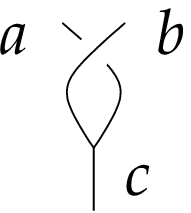}{1.3cm}=\d(c; a,b) \ppsdraw{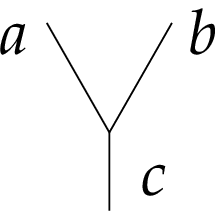}{1.3cm}
\qquad
\ppsdraw{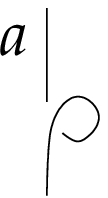}{1.3cm}=\mu(a) \ppsdraw{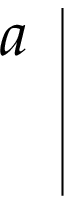}{1.3cm}
$$
Here, $\la i,j,k \ra$ is defined below and 
\begin{eqnarray*}
\la k \ra &=& (-1)^k[k+1]=(-1)^k\frac{A^{2k+2}-A^{-2k-2}}{A^2-A^{-2}}
\\
\mu(k) &=& (-1)^k A^{k^2+2k} \\
\d(c;a,b) &=& (-1)^k A^{ij-k(i+j+k+2)}.
\end{eqnarray*}

With the use of fusion one can reduce the computation of knotted trivalent
graphs to the values of a standard trihedron and tetrahedron.

Following standard skein theory conventions, a triple $(a,b,c)$ of
nonnegative integers is called {\em admissible} if $a+b+c$ is even
and $|b-c| \leq a \leq b+c$. Given an admissible triple, let
$i=(b+c-a)/2$, $j=(a+c-b)/2$ and $k=(a+b-c)/2$ denote the corresponding
{\em internal colors}.
Of importance are the values of a colored trihedron 
$$
\psdraw{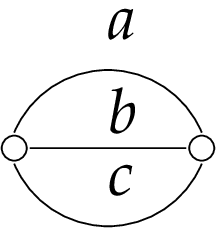}{0.5in}=\la a,b,c \ra
$$
and tetrahedron
$$
\psdraw{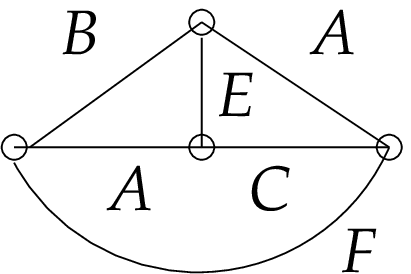}{1in}=
\left\la \begin{matrix} A & B & E \\ D & C & F \end{matrix}\right\ra.
$$
Using the notation of Masbaum-Vogel \cite{MV} and also Kauffman-Lins
\cite{KL}, the trihedron and tetrahedron coefficients are given by:
\begin{eqnarray*}
\la a,b,c \ra &=& (-1)^{i+j+k}\frac{[i+j+k+1]![i+j]![i+k]![j+k]!}{[a]![b]![c]!}
\\
\left\la \begin{matrix} A & B & E \\ D & C & F \end{matrix}\right\ra
& = & \frac{\prod_{i=1}^3\prod_{j=1}^4 [b_i-a_j]!}{
[A]![B]![C]![D]![E]![F]!} 
\left( \begin{matrix} 
a_1 &     & a_2 &     & a_3 &      & a_4 \\
   & b_1 &     & b_2 &     & b_3  &
\end{matrix}
\right),
\end{eqnarray*}
where
\begin{itemize}
\item
we assume that the triples $(A,B,E), (B,D,F), (E,D,C)$ and $(A,C,F)$ are 
admissible,
\item
$\S=A+B+C+D+E+F$, 
\begin{align*}
a_1 &=(A+B+E)/2 & b_1 &=(\S-A-D)/2 \\
a_2 &=(B+D+F)/2 & b_2 &=(\S-E-F)/2 \\
a_3 &=(C+D+E)/2 & b_3 &=(\S-B-C)/2 \\
a_4 &=(A+C+F)/2, &     &
\end{align*}\end{itemize}
$\phantom{999}\bullet$\qua and 
$$
\left( \begin{matrix} 
a_1 &     & a_2 &     & a_3 &      & a_4 \\
   & b_1 &     & b_2 &     & b_3  &
\end{matrix}
\right)=\sum_{\text{max}(a_j) \leq \z \leq \text{min}(b_i)}
\frac{(-1)^{\z}[\z+1]!}{\prod_{i=1}^3[b_i-\z]!\prod_{j=1}^4 [\z-a_j]!}.
$$

\subsection{Proof of theorem \ref{thm.CJ}}
\label{sub.proof}

Recall that $K_0$ is the closure of the braid $b^2(ab)^8=ba^{-1}(ab)^9=ba^{-1}
c^3$ where $c=(ab)^3$ represents a full twist.
$c$ can be obtained by first giving a full twist on the first two strands,
and then a full twist of the third strand around the first two.
Thus, $K_0$ is obtained from the closure of the 
following figure:
$$
\ppsdraw{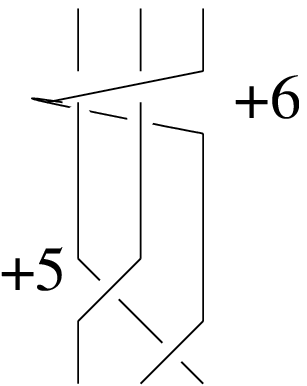}{2cm}
$$
We want to compute the Kauffman bracket $\la K_0,n \ra$ colored by $n$, using
the zero framing of $K$; this differs from the blackboard framing by $+18$.
Now, let us fuse the first two strands, and undo the $+5$ half twists.
We obtain that
$$
\la K_0,n\ra=\mu(n)^{-18}
\sum_k \frac{\la k \ra}{\la n,n,k\ra} \ppsdraw{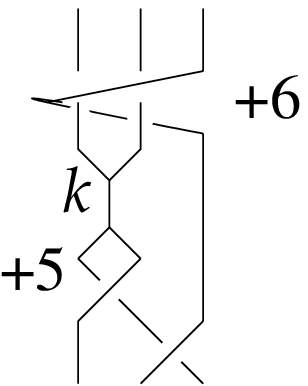}{2cm}
=
\sum_k \frac{\la k \ra}{\la n,n,k\ra} \delta(k;n,n)^5 \ppsdraw{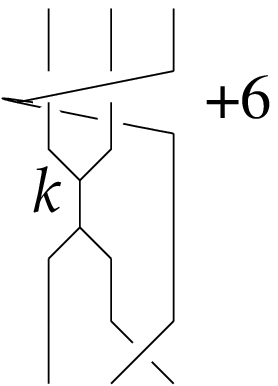}{2cm}
$$
where the color of any noncolored edges is $n$.
Let us isotope the fused strand above the $+6$ half twists, fuse again, and
then undo the $+6$ half-twists. We obtain that
$$
\la K_0,n\ra=\mu(n)^{-18}
\sum_{k,l} \frac{\la k \ra}{\la n,n,k\ra}
\frac{\la l \ra}{\la l,k,n\ra} \delta(k;n,n)^5
\delta(l;k,n)^6
\ppsdraw{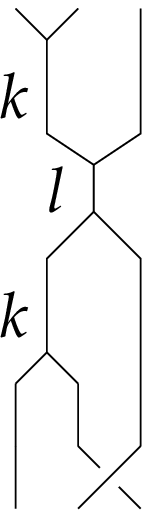}{3cm}
$$
Now, let us close up to a tetrahedron. Taking into account two half-twists,
we obtain that
$$
\la K_0,n\ra=\mu(n)^{-18}
\sum_{k,l} \frac{\la k \ra}{\la n,n,k\ra}
\frac{\la l \ra}{\la l,k,n\ra} \delta(k;n,n)^4
\delta(l;k,n)^7
\left\la \begin{matrix} n & n & k \\ n & l & k \end{matrix}\right\ra
$$
Now, let us define the renormalized trihedron and tetrahedron
coefficients $\la \cdot \ra'$ by:
\begin{eqnarray*}
\left\la \begin{matrix} A & B & E \\ D & C & F \end{matrix}\right\ra &=&
\frac{1}{[A]![B]![C]![D]![E]![F]!} 
\left\la \begin{matrix} A & B & E \\ D & C & F \end{matrix}\right\ra' \\
\la a,b,c \ra &=& \frac{1}{[a]![b]![c]!} \la a,b,c \ra'
\end{eqnarray*}
Since the
colors on the edges of the tetrahedron cancel the colors of the edges
of the trihedron, we obtain
$$
\la K_0,n\ra=\mu(n)^{-18}
\sum_{k,l} \frac{\la k \ra}{\la n,n,k\ra'}
\frac{\la l \ra}{\la l,k,n\ra'} \delta(k;n,n)^4
\delta(l;k,n)^7
\left\la \begin{matrix} n & n & k \\ n & l & k \end{matrix}\right\ra'
$$
In the above expression,
$
\mu(n)^{-18} \delta(k;n,n)^4
\delta(l;k,n)^7$ is a monomial of $A$. 

The $a_i$ and $b_j$'s for the tetrahedron in question are given by:
\begin{eqnarray*}
a_1=a_3 &=& (2n+k)/2 \\
a_2=a_4 &=& (n+k+l)/2 \\
b_1 &=& (n+2k+l)/2 \\
b_2 &=& (3n+l)/2 \\
b_3 &=& n+k.
\end{eqnarray*}
Since $a_i-b_j$ are the internal colors at the vertices of the tetrahedron,
it follows that $6$ quantum factorials cancel. In other words, we have:
$$
\frac{\la k \ra}{\la n,n,k\ra'}
\frac{\la l \ra}{\la l,k,n\ra'} \prod_{i=1}^3\prod_{j=1}^4 [b_i-a_j]!
=\prod_{i=1}^3 [b_i-a_1]! \prod_{i=1}^3 [b_i-a_3]!
$$
Combining the remaining $6$ quantum factorials, and using the fact that
$$
J_K(n+1)=\frac{(-1)^n}{[n+1]}\la K,n\ra,
$$ 
the result follows.
\qed

Remark \ref{rem.qbinom} 
follows from the fact that the quantum binomial coefficients
are Laurent polynomials.

\section{The classical topology and geometry of $K_0$}
\label{sec.classical}

In this largely independent section, we give several facts about the classical
geometry and topology of $K_0$.

\subsection{The topology of $K_0$}
\label{sub.top}

Torus knots can be embedded in closed surfaces of genus $1$.
Likewise, twisted torus knots can be embedded in closed surfaces of genus $2$.
$K_0$ is a $(1,1)$ {\em knot}, that is a special kind of a {\em tunnel number 
$1$ knot}; see for example \cite{GMM} and \cite{MSY}.
Consequently, its fundamental group has {\em rank} $2$, i.e., 
is a $2$-generator and
$1$-relator group, where the meridian is one of the generators.
Note that $2$-bridge knots are also $(1,1)$ knots; thus from the point
of view of presentation of the fundamental group, the knot $K$ looks similar
to a $2$-bridge knot, although it is not one.
Dean studies which twisted torus knots are 
$(1,1)$ knots, \cite{De}.

SnapPea reveals that the fundamental group of $K_0$ has the following 
presentation:

$$
\pi_1(S^3-K_0)=\la a,b | 
aBabaBabbaBaBabbaBaBabb \ra
$$

where the canonical meridian-longitude pair $(m,l)$ is given by
$$
(m,l)= (a, ABBAbABB),
$$ 
where $A=a^{-1}$ and $B=b^{-1}$.

%

The Alexander polynomial of $K_0$ is:
$$
t^{-8}-t^{-7}+t^{-5}-t^{-4}+t^{-2}-t^{-1}+1-t+t^2-t^4+t^5-t^7+t^8
$$
and the sum of the absolute value of its coefficients is $13$. 
By \cite{GMM}, and a private conversation with J. Rasmussen, $13$ is also the 
rank of the {\em Heegaard Floer Homology} of Oszvath-Szabo; \cite{OS}.

\subsection{The geometry of $K_0$}
\label{sub.hypgeom}

SnapPea numerically computes the volume of $K_0$ to be:
$$
\vol(K_0)=  3.474247 \ldots .
$$

Of importance to us is the $A$-{\em polynomial} of $K_0$, which parametrizes
the complex $1$-dimensional part of the $\SL_2(\BC)$ character variety, 
restricted to the boundary of the knot complement. The $A$-polynomial of
a knot was introduced by \cite{CCGLS}, and is important in 
\begin{itemize}
\item[(a)] the study of
deformations of the complete hyperbolic structure, 
\item[(b)] the study of slopes of essential surfaces in the knot
complement with nonempty boundary, and 
\item[(c)] in the shape of linear recursion relations of the colored Jones
function, \cite{Ga}.
\end{itemize}

Since $K_0$ has rank $2$, with a bit of work one can compute the $A$-polynomial
of $K_0$, as was explained in \cite{CCGLS}. The method of \cite[Section 7]{CCGLS}
assigns $2$ by $2$ matrices to the generators of the fundamental group,
and uses resultants to compute the polynomial relation satisfied by
the eigenvalues of the meridian and the longitude.

Alternatively, SnapPea provides an ideal triangulation of the knot complement
of $K_0$, and one can write down the corresponding {\em gluing equations}
and use elimination (via Groebner basis) to compute 
the polynomial relation satisifed by the eigenvalues of the meridian and 
the longitude. Boyd has developed a few more
methods for computing the $A$-polynomial once its expected degree is known.
Either method gives that
$A(l,m)=(l-1)B(l,m^2)$ where $B(l,m)$ is an irreducible 
polynomial of $(l,m)$ of degree $10 \times 136$. Explicitly, we have:
\begin{align*}
&B(l,m)=\\
&\quad\  l^{10} + l^9 m^{13} + l^9 m^{14} - l^8 m^{25} + 4 l^8 m^{26} - 
    8 l^8 m^{27} + 3 l^8 m^{28} - l^8 m^{29} - 3 l^7 m^{40} \\ & - 
    4 l^7 m^{41} - l^7 m^{42} + 2 l^6 m^{55} - 2 l^5 m^{66} + 
    11 l^5 m^{67} - 6 l^5 m^{68} + 11 l^5 m^{69} - 2 l^5 m^{70} \\ & + 
    2 l^4 m^{81} - l^3 m^{94} - 4 l^3 m^{95} - 3 l^3 m^{96} - 
    l^2 m^{107} + 3 l^2 m^{108} - 8 l^2 m^{109} + 4 l^2 m^{110} \\ & - 
    l^2 m^{111} + l m^{122} + l m^{123} + m^{136}.
\end{align*}
We are interested in the $10$ roots $\l_k(t)$ (so-called {\em eigenvalues})
of the equation 
$$
B(\l,e^{i t})=0.
$$
We have that:
$$
B(l,1)=(l+1)^6 (l-1)^4.
$$ 
Since $B$ is reciprocal with real coefficients, it follows that if $\l_k(t)$
is an eigenvalue, so is the complex conjugate $\bar{\l}_k(t)$ and the inverse
$1/\l_k(t)$. Thus, the set of the eigenvalues has $4$-fold symmetry. 
It turns out that $6$ eigenvalues $\l_k(t)$  for $t=5,\dots,10$
have magnitude $1$ for all $t \in [0,2\pi]$. The remaining four eigenvalues
are 
$$
\l_1(t), \qquad \l_2(t)=\l_1(2\pi-t), \qquad \l_3(t)=\bar{\l}_1(t),
\qquad \l_4(t)=\bar{\l}_1(2\pi-t),
$$
where $\l_1$ is chosen so that the slope of its magnitude at $0$ is bigger
than the slope of the magnitude of $\l_2$ at $0$. 

The plot of the set $\{ \log|\l_k(t)|\, | k=1,\dots,10 \}$
versus $t$ is:
$$
\psdraw{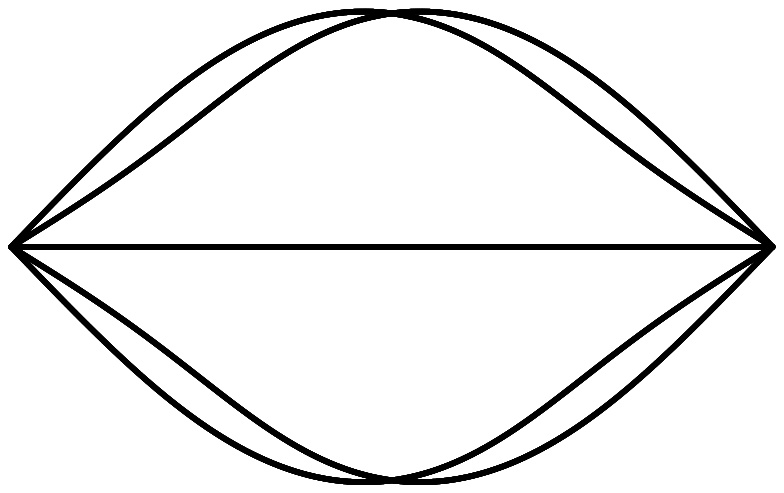}{2in}
$$
Collision of eigenvalues occurs only at $t=0$, $t=\pi$ and $t=2\pi$,
as follows from computing the roots of the discriminant 
$\text{Disc}_lB(l,m)$ which lie on the unit circle.

Using the formula for the variation of the volume in terms of the 
$A$-polynomial, it follows that 
$$
\int_0^{\pi} \log|\l_1(t)| \,dt = V_1 + V_3
\qquad
\int_0^{\pi} \log|\l_2(t)| \,dt = V_2 - V_3
$$
where
$$
V_1=3.474247 \ldots,  \qquad V_2=1.962737 \ldots, \qquad V_3=0.490684 \ldots
$$
Here, $V_i=\text{vol}(\rho_i)$ for $i=1,2,3$ where
\begin{itemize}
\item
$\rho_1$ is the discrete faithful representation, passing through 
$(m,l)=(1,-1)$. 
\item
$\rho_2$ is another solution, passing though $(m,l)=(1,1)$,
\item
$\rho_3$ is a discrete faithful of the hyperbolic Dehn filling $K(2,0)$.
\end{itemize}

Using Snap, Boyd informs us that the 
{\em invariant trace fields} $E_i$ of $\rho_i$ are of type
\begin{itemize}
\item
$6, \,\, [4,1], \,\, -753079$ for $\rho_1$, (that is, of degree $6$, with
$4$ real embeddings, $1$ complex embedding and conductor $-753079$),
\item
$4, \,\, [2,1], \,\, -283$ for $\rho_2$ and $\rho_3$.
\end{itemize}

It is a coincidence that $E_2=E_3$, which implies that $V_3=V_2/4$.

The entropy of a knot was defined in \cite{GG}. Using the above information,
it follows that the possible values of the entropy (at $t=2\pi$)
are given by the areas of the following curves:
$$
\psdraw{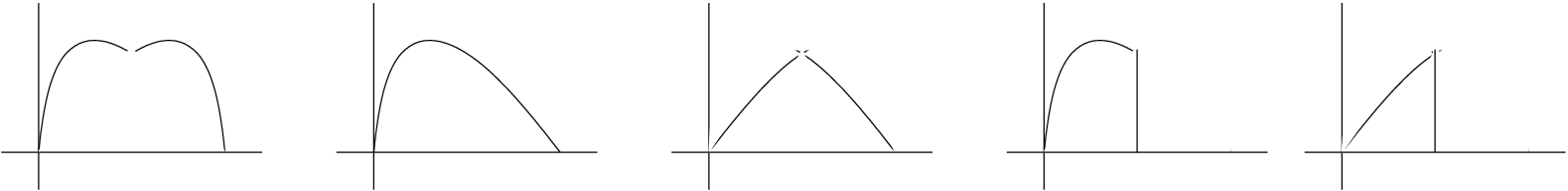}{.9\hsize}
$$
(or their reflection about $t=\pi$). These values are, respectively:
\begin{eqnarray*}
2(V_1+V_3) &=& 7.929864 \ldots \\ 
V_1+V_2 &=& 5.436985 \ldots \\ 
2(V_2-V_3)&=& 2.94411 \ldots \\ 
V_1+V_3&=& 3.964932 \ldots \\ 
V_2-V_3 &=& 1.472053 \ldots .
\end{eqnarray*}
Notice that none of these values equals to the volume $V_1$ of $K_0$.

\section{Towards the AJ conjecture for the knot $K_0$}
\label{sec.AJ}

It was proven in \cite{GL} that the sequence $J_K(n)$ for $n=0,1,2,\dots$
of Laurent polynomials is $q$-holonomic, that is it satisfies a recursion
relation.

Actually, it was explained in \cite{GL} that the colored Jones function
of any knot is given as a multisum of a $q$-proper hypergeometric function.
Using the general theory of Wilf-Zeilberger, it follows that such multisums
always satisfy recursion relations.

In general, the multisums in \cite{GL} use as many summation variables as
the number of the crossings of the knot, which makes things impractical
to compute. In our case, we may use Theorem \ref{thm.CJ} which is a triple
sum only.

What can one say about a recursion relation for $J_{K_0}$? It was conjectured
in \cite{Ga} (as part of the AJ Conjecture) that $J_{K_0}$ satisfies a degree
$11$ recursion relation of the form:
\begin{equation}
\label{eq:AJ}
\sum_{j=0}^{11} a_j(q^n,q) J_{K_0}(n+j)=0
\end{equation}
where
\begin{itemize}
\item
$a_j(u,v) \in \BQ(u,v)$ are rational functions with rational coefficients,
\item
the recursion relation is a $q$-{\em deformation} 
of the $A$-polynomial. That is,
$$
\sum_{j=0}^{11} a_j(M^2,1) L^j=A(L,M).
$$
\item
$\sum_{j=0}^{11} a_j(Q,q)E^j$ is {\em reciprocal} in the sense that
$$
a_j(Q,q)=\e_j a_{11-j}(q^{-11}Q^{-1},q)
$$
for all $j=0,\dots,11$ where $\e_j =\pm 1$. 
\end{itemize}
These properties cut down a guess for a recursion relation considerably.

Using Theorem \ref{thm.CJ}, we have computed explicitly $J_{K_0}(n)$ for
$n=1,\dots,19$. It is a challenging question to guess a formula for the
so-called noncommutative $A$ polynomial of the knot $K_0$.




\subsection{Acknowledgement}

Many people contributed to the ideas and the computations of this paper.
We wish to thank D. Boyd, A. Champanerkar and N. Dunfield 
(for counseling on the $A$-polynomial), 
T. Morley (for {\tt Mathematica} tutoring), M. Thistlethwaite (for
the symmetry of $K$) and 
TTQ. Le (for discussions on a triple sum formula for the colored Jones
function).

\appendix

\section{The numerical computation of $\VC(n)$}

\subsection{A recursive computation of $\VC(n)$}
In this appendix, we discuss the numerical evaluation 
of the colored Jones polynomial $J(n+1):=J_{K_0}(n+1)$ at $q=q_0$, where 
$q_0=\exp(i\frac{2\pi}{n+1})$. 
Let us begin by writing the formula in Theorem \ref{thm.CJ} in a numerically 
more efficient form
\begin{align}
&J(n+1)=\nonumber\\
 &\quad \frac{1}{[n+1]} \sum_{k=0,2}^{2n}\sum_{l=|n-k|,2}^{n+k}
q^{-\frac{3}{8}(2k+k^2)+\frac{7}{8}(2l+l^2)-\frac{51}{8}(2n+n^2)}
\frac{[k+1][l+1]\left[\frac{k}{2}\right]!}{[\frac{2n+k}{2}+1,\frac{k}{2}]!} 
Z \,,\nonumber\\
\intertext{with\ \ $Z=$}
 & \sum_z (-1)^{\frac{k}{2}+z}  
\qb {\frac{k+l-n}{2}}{\frac{n+2k+l}{2}-z}
\qb {\frac{n+l-k}{2}}{\frac{3n+l}{2}-z}
\qb {\frac{n+k-l}{2}}{\frac{2n+2k}{2}-z}  
\left[\frac{2n-k}{2},z-\frac{n+k+l}{2}\right]! \nonumber\\
&\qquad\left[z+1,\frac{n+k+l}{2}+1\right]! 
\,,\label{eq:jonesp}
\end{align}
so that all the factors in $Z$ have $z$ dependence. 

If we define $B=q^{1/2}$, then $Z$ in Equation \eqref{eq:jonesp} 
is an explicit rational function of $B$. We know that $J(n+1)$ is a finite 
sum with integer coefficients of integer
powers of $q$, so it is also a finite sum of even powers of $B$. 
From Theorem \ref{thm.CJ} and Remark \ref{rem.qbinom}, we know that as a 
function of $B$, $J(n+1)(B)$ has an apparent second order pole at
$B=B_0=\exp(i\frac{\pi}{n+1})$. To verify the Volume Conjecture, we need to 
evaluate the Jones polynomial of the form:
\begin{equation}
J(n+1)(B_0)=\frac{f(B)}{(B-B_0)^2}|_{B= B_0}\,,
\end{equation}
where $f(B)=g(B)(B-B_0)^2$ with $g(B)$ an analytic function at 
$B=B_0$. Apparently, $J(n+1)(B_0)=g(B_0)$. On the other hand
\begin{equation}
g(B_0)=\frac{1}{2}\frac{d^2}{dB^2}f(B)|_{B=B_0}
\,.
\end{equation}
The next Equation will be the basis of our numerical
computation:
\begin{equation}
\label{eq:compute}
J(n+1)(B_0)=\frac{1}{2}\frac{d^2}{dB^2}f(B)|_{B=B_0}.
\end{equation}
In other words, we have to separate the factor $(B-B_0)^2$ in the denominator 
in Equation \eqref{eq:jonesp}, differentiate the rest part twice and evaluate 
at $B=B_0$. Considering the three summations and the complicated structure of
the summand, we would expect that for large $n$, the computation load is
huge. This is true if we employ symbolic manipulations using {\em Mathematica}
or {\em Maple}, but numerically, we can do it much more efficiently if
recursion relation is used both for the quantum numbers and their 
derivatives.

We have to evaluate the {\em quantum integer}, the {\em quantum factorial} 
and the {\em quantum binomial coefficients} etc. defined just before
Theorem \ref{thm.CJ}. Then 
the {\em quantum integer} is
\begin{equation}
[n]=\frac{B^n-B^{-n}}{B-B^{-1}}
\,,
\end{equation}
which has  first derivative
\begin{equation}
[n]^{\prime}=\frac{n(B^{n-1}+B^{-n-1})}{B-B^{-1}}-\frac{B^n-B^{-n}}
{(B-B^{-1})^2}(1+B^{-2})
\,,
\end{equation}
and second derivative
\begin{align}
[n]^{\prime \prime}=&\frac{n((n-1)B^{n-2}-(n+1)B^{-n-2})}{B-B^{-1}}
-\frac{2n(B^{n-1}+B^{-n-1})}{(B-B^{-1})^2}(1+B^{-2})\nonumber\\
&+\frac{2(B^n-B^{-n})}{(B-B^{-1})^3}(1+3B^{-2})
\,.
\end{align}
In numerical calculations, three vectors may be used to store these values.
 
The {\em quantum factorial} satisfies the recursion relation
\begin{equation}
[n]!=[n][n-1]!
\,,
\end{equation}
which induces recursion relations for its first derivative 
\begin{equation}
[n]!^{\prime}=[n]^{\prime}[n-1]!+[n][n-1]!^{\prime}
\end{equation}
and second derivative
\begin{equation}
[n]!^{\prime \prime}=[n]^{\prime \prime}[n-1]!+2[n]^{\prime}[n-1]!^{\prime}
+[n][n-1]!^{\prime \prime}
\,.
\end{equation}
Starting from $[1]!=1$, these factorials can be calculated simultaneously 
with the {\em quantum 
integers} and also stored in three vectors.

Similarly, for {\em quantum binomial coefficients}, the recursion relations 
are
\begin{equation}
\qb n k=B^{-k} \qb {n-1} k + B^{n-k} \qb {n-1}{k-1}
\,,
\end{equation}
\begin{equation}
\qb n k=B^{-k} \qb {n-1} k ^{\prime}+ B^{n-k} \qb {n-1}{k-1} ^{\prime}
-kB^{-k-1} \qb {n-1} k + (n-k)B^{n-k-1} \qb {n-1}{k-1}
\,,
\end{equation}
and\ \ $\qb n k=$ 
\begin{align*}
& B^{-k} \qb {n-1} k ^{\prime \prime}+ B^{n-k} \qb {n-1}{k-1} 
^{\prime \prime}-kB^{-k-1} \qb {n-1} k ^{\prime}+ (n-k)B^{n-k-1} 
\qb {n-1}{k-1} ^{\prime}  \\
&+ k(k+1)B^{-k-2} \qb {n-1} k + (n-k)(n-k-1)B^{n-k-2} \qb {n-1}{k-1}
\,,
\end{align*}
with the starting condition $\qb n n =1$.
 
The recursion relations related to $[n,k]!$ are given by
\begin{equation}
[n,k]!=[n][n-1,k]!
\,,
\end{equation}
\begin{equation}
[n,k]!^{\prime}=[n]^{\prime}[n-1,k]!+[n][n-1,k]!^{\prime}
\,,
\end{equation}
and
\begin{equation}
[n,k]!^{\prime}=[n]^{\prime \prime}[n-1,k]!+[n][n-1,k]!^{\prime \prime}
+[n]^{\prime}[n-1,k]!^{\prime}
\,,
\end{equation}
with starting condition $[k,k]!=1$. Notice that for the {\em quantum 
binomial coefficients} $\qb n k$ and $[n,k]!$, $ n\geq k $ has to be
satisfied. So, their numerical values can be put together to form
a square matrix.

We need to take special care of the denominator $D$ in each summand. It has
the form  
\begin{equation}
D=[\frac{2n+k}{2}+1,\frac{k}{2}]![n+1]=[\frac{2n+k}{2}+1,n+1]!
[n,\frac{k}{2}]!P
\,,
\end{equation}
where $P=[n+1]^2$.
Note that $D$ is factorized to three factors. The first two with their 
first and second derivatives are readily evaluated by previous calculations
and do not have $q=q_0$
as zeros. All the zero factors are contained in the third factor P. We have
to extract the factor $B-B_0$ and evaluate the rest.
we may write
\begin{equation}
[n+1]=\frac{B^{n+1}-B^{-n-1}}{B-B^{-1}}=TS(B-B_0)
\,,
\end{equation}
where
\[
T=\frac{B^{-n-1}}{B-B^{-1}}\,, \qquad S=B^{2n+1}+B^{2n}B_0+\cdots+B_0^{2n+1}
\,.
\]
It is straight forward to compute that
\begin{eqnarray*}
S(B_0)&=&2(n+1)B_0^{-1} \\
S^{\prime}(B_0)&=&(n+1)(2n+1)B_0^{-2} \\
S^{\prime \prime}(B_0)&=&\frac{4}{3}n(n+1)(2n+1)B_0^{-3} 
\,;
\end{eqnarray*}
and
\begin{eqnarray*}
T^{\prime}(B)&=&\frac{-(n+1)B^{-n-2}}{B-B^{-1}}-\frac{B^{-n-1}}{(B-B^{-1})^2}
(1+B^{-2}) \\
T^{\prime\prime}(B)&=&\frac{(n+1)(n+2)B^{-n-3}}{B-B^{-1}}+\frac{2(n+1)B^{-n-2}}
{(B-B^{-1})^2}(1+B^{-2})\\&&+\frac{2B^{-n-1}}{(B-B^{-1})^3}(1+3B^{-2}) 
\,.
\end{eqnarray*}
Now, we are in a position to compute the first and second 
derivatives of the summand of Equation 
\eqref{eq:jonesp}, by invoking the values which we have.
Finally, $g(B_0)$ is computed by 
multiplication and addition of the terms.

\subsection{The multiprecision algorithm}

Using the above scheme, we can go up to $n\sim 200$ with double precision in
Fortran. Beyond that, the calculation overflows. The sequence $\VC(n)$ 
decreases until $n\sim60$, where we reach a 
minimum, then it starts to increase. Of course, this is only
an illusion due to the numerical round-off error. We checked the magnitude of
the terms in the sum, it could go beyond that of the final result, 
which suggests that 
there could be {\em significant cancellations} among terms and the precision 
of our calculation is below the required one. We have to use multi-precision
algorithm.

We take advantage of the multiprecision program by Bailey~\cite{Ba}
(available at \url{http://crd.lbl.gov/~dhbailey/mpdist/}). 
We used 
MPFUN90 (Fortran-90 arbitrary precision package) to do our calculation.
Using $40$ digits, the series reaches a minimum at $n=131$. With $80$
digits, the summation can be easily done if $n\sim250$. The results shown
in section \ref{sec.numerical} are obtained with $80$ digits. we also 
tried $200$
digits and calculated values up to $n=550$. Beyond that, the computer ran out
of memory. We believe that with larger memory, it is easy for us to go beyond
$n=1000$. With the increase of digits, we can go higher and higher values of
$n$ but with the price that the computation is done slower and slower. \\

All the calculations that we have done suggest that the results in section
\ref{sec.numerical} is true with high credibility. We may use more digits to
calculate the $\VC(n)$ to ever-higher precision if necessity arises.

\section{The colored Jones function of $K_0$}

Using Theorem \ref{thm.CJ} we can compute the colored Jones function
$J(n):=J_{K_0}(n)$ for $n=1,\dots,19$. 
The polynomials become large soon, and we list here the first few of them.
\begin{align*}
&J(1) = 1 \\
&J(2) =
-\frac{1}{{q^{19}}}+\frac{1}{{q^{18}}}-\frac{1}{{q^{17}}}
+\frac{1}{{q^{10}}}+\frac{1}{{q^8}} 
\\
&J(3)=\\
&\frac{1}{{q^{53}}}\big(1-q+{q^2}-2 {q^4}+2 {q^5}-2 {q^7}+{q^8}+{q^9}
-{q^{10}}+{q^{12}}- {q^{13}}-{q^{14}}+{q^{15}}-{q^{17}}
\\ & 
+{q^{19}}
-{q^{20}}-{q^{21}}+{q^{22}}-{q^{24}}+{q^{25}}- {q^{27}}+{q^{28}}
-{q^{30}}+{q^{31}}+{q^{34}}+{q^{37}}\big)
\end{align*}
\begin{align*}
&J(4) =\\&
\frac{1}{{q^{100}}}\big(-1+2  {q^2}-4  {q^4}+{q^5}+4  {q^6}-5  {q^8}+6  
{q^{10}}-5  {q^{12}}+
4  {q^{14}}+{q^{15}}-4  {q^{16}}+3  {q^{18}}  \\
& +{q^{19}}-3  {q^{20}}-{q^{21}}+2  {q^{22}}+  
{q^{23}}-{q^{24}}-2  {q^{25}}+{q^{26}}+{q^{27}}-{q^{29}}
-{q^{30}}+{q^{31}}+{q^{32}} \\ & -  
{q^{34}}+{q^{36}}-{q^{38}}+{q^{40}}-{q^{41}}-{q^{42}}
+2  {q^{44}}-{q^{45}}-{q^{46}}-  {q^{47}}+2  {q^{48}}-{q^{50}}-{q^{51}}
\\ & +{q^{52}}-{q^{54}}+{q^{56}}-{q^{58}}+{q^{60}}-  
{q^{62}}+{q^{64}}-{q^{66}}+{q^{68}}+{q^{72}}+  {q^{76}}\big)
\end{align*}
\begin{align*}
&J(5) =\\&
\frac{1}{{q^{160}}}\big(2-q-{q^2}-2 {q^3}+5 {q^5}+{q^6}-3 {q^7}
-6 {q^8}+7 {q^{10}}+   
 4 {q^{11}}-3 {q^{12}}-8 {q^{13}}-2 {q^{14}} 
\\ & +7 {q^{15}}+6 {q^{16}}-2 {q^{17}}-8
{q^{18}}-   
 2 {q^{19}}+6 {q^{20}}+4 {q^{21}}-2 {q^{22}}-6 {q^{23}}-2 {q^{24}}+6 {q^{25}}
\\ & 
+   
 3 {q^{26}}-{q^{27}}-5 {q^{28}}-2 {q^{29}}+5 {q^{30}}+2 {q^{31}}
-3 {q^{33}}-2 {q^{34}}+  
 3 {q^{35}}+{q^{36}}+{q^{37}}-{q^{38}}
\\ &
-2 {q^{39}}+{q^{40}}+{q^{42}}+{q^{43}}-{q^{44}}-   
 {q^{46}}+{q^{48}}+{q^{49}}-{q^{51}}-{q^{52}}+{q^{54}}+{q^{56}}-{q^{57}}
\\ &
-{q^{58}}-{q^{60}}+   
 2 {q^{61}}-3 {q^{65}}+2 {q^{66}}+{q^{68}}+{q^{69}}-3 {q^{70}}+2 {q^{71}}
-{q^{72}}+   
 {q^{73}}+{q^{74}}
\\ &
-3 {q^{75}}+2 {q^{76}}-{q^{77}}+{q^{78}}+{q^{79}}
-3 {q^{80}}+{q^{81}}-   
 {q^{82}}+{q^{83}}+{q^{84}}-2 {q^{85}}+{q^{86}}
\\ &
-{q^{87}}+{q^{88}}
-2 {q^{90}}+{q^{91}}+   
 {q^{93}}-2 {q^{95}}+{q^{98}}-{q^{100}}+{q^{103}}-{q^{105}}+{q^{108}}
-{q^{110}}
\\ &
+   
 {q^{113}}-{q^{115}}+{q^{118}}+{q^{123}}+{q^{128}}\big)
\end{align*}
\begin{align*}
&J(6) =\\&
\frac{1}{{q^{235}}}\big(-1+{q^3}+{q^4}-{q^7}-2 {q^8}-{q^9}+5 {q^{11}}
+4 {q^{12}}-{q^{13}}-   
 5 {q^{14}}-6 {q^{15}}-{q^{16}}
\\ &
+6 {q^{17}}+9 {q^{18}}+2 {q^{19}}-6 {q^{20}}-9 {q^{21}}-
 4 {q^{22}}+4 {q^{23}}+9 {q^{24}}+6 {q^{25}}-5 {q^{26}}
-8 {q^{27}}
\\ &
 -4 {q^{28}}+   
 3 {q^{29}}+8 {q^{30}}+3 {q^{31}}-4 {q^{32}}-7 {q^{33}}-3 {q^{34}}+4 {q^{35}}+7
{q^{36}}+   
 2 {q^{37}}-3 {q^{38}}
\\ &
-6 {q^{39}}-3 {q^{40}}+3 {q^{41}}+5 {q^{42}}+2 {q^{43}}-2
{q^{44}}-   
 3 {q^{45}}-3 {q^{46}}+{q^{47}}+3 {q^{48}}+2 {q^{49}}
\\ &
-{q^{51}}
-2 {q^{52}}-3 {q^{53}}+   
 {q^{54}}+2 {q^{55}}+2 {q^{56}}+{q^{57}}-{q^{58}}-4 {q^{59}}
-2 {q^{60}}+{q^{61}}+   
 3 {q^{62}}
\\ &
+3 {q^{63}}+{q^{64}}-3 {q^{65}}-3 {q^{66}}-2 {q^{67}}
+2 {q^{68}}+4 {q^{69}}+
 2 {q^{70}}-3 {q^{72}}-3 {q^{73}}+2 {q^{75}}
\\ &
+2 {q^{76}}+2 {q^{77}}
-{q^{78}}-2 {q^{79}}- {q^{81}}+{q^{82}}+{q^{83}}+{q^{86}}-{q^{87}}
-{q^{90}}+{q^{92}}+{q^{94}}
\\ &
-{q^{96}}-{q^{99}}+   
 {q^{100}}-{q^{102}}+2 {q^{103}}-{q^{105}}-{q^{107}}-2 {q^{108}}
+2 {q^{109}}+   {q^{110}}+{q^{112}}
\\ &
-{q^{113}}-2 {q^{114}}+{q^{116}}
+{q^{118}}-{q^{120}}+{q^{122}}- {q^{126}}+{q^{128}}-{q^{132}}+{q^{134}}
-{q^{136}}
\\ &
-{q^{138}}+{q^{140}}+{q^{141}}- {q^{142}}-{q^{144}}+{q^{147}}
-{q^{148}}+{q^{153}}-{q^{154}}-{q^{155}}+{q^{159}}
\\ &
- {q^{161}}
+{q^{165}}-{q^{167}}+{q^{171}}-{q^{173}}+{q^{177}}-{q^{179}}+{q^{183}}+   
 {q^{189}}+{q^{195}}\big)
\end{align*}
\begin{align*}
&J(7) =\\&
\frac{1}{{q^{323}}}\big(1+{q^2}-3 {q^3}+3 {q^7}-{q^8}-5 {q^{10}}
+3 {q^{11}}+2 {q^{12}}+{q^{13}}+ 3 {q^{14}}-6 {q^{15}}-5 {q^{16}}
\\ &
-5 {q^{17}}+9 {q^{18}}+8 {q^{19}}+3 {q^{20}}+3 {q^{21}}-   
 11 {q^{22}}-13 {q^{23}}-8 {q^{24}}+12 {q^{25}}+12 {q^{26}}
\\&+9 {q^{27}}
+7 {q^{28}}-  12 {q^{29}}-17 {q^{30}}-13 {q^{31}}+9 {q^{32}}+11 {q^{33}}
+10 {q^{34}}+12 {q^{35}}- 9 {q^{36}}
\\ &
-16 {q^{37}}-13 {q^{38}}
+8 {q^{39}}
+10 {q^{40}}+8 {q^{41}}+9 {q^{42}}- 9 {q^{43}}-14 {q^{44}}-11 {q^{45}}
+9 {q^{46}}
\\ &
+10 {q^{47}}+8 {q^{48}}+7 {q^{49}}
- 9 {q^{50}}-12 {q^{51}}
-11 {q^{52}}+7 {q^{53}}+8 {q^{54}}+7 {q^{55}}+6 {q^{56}}
\\ &
- 6 {q^{57}}
-8 {q^{58}}-11 {q^{59}}+4 {q^{60}}
+6 {q^{61}}+6 {q^{62}}+7 {q^{63}}
- 3 {q^{64}}-5 {q^{65}}-11 {q^{66}}-{q^{67}}
\\ &
+3 {q^{68}}+6 {q^{69}} +8 {q^{70}}+{q^{71}}- 2 {q^{72}}
-11 {q^{73}}-4 {q^{74}}-{q^{75}}
+5 {q^{76}}+8 {q^{77}}+4 {q^{78}}
\\ &
+2 {q^{79}}- 9 {q^{80}}-6 {q^{81}}
-4 {q^{82}}+2 {q^{83}}+6 {q^{84}}
+5 {q^{85}}+5 {q^{86}}- 4 {q^{87}}
-5 {q^{88}}-6 {q^{89}}
\\ &
-{q^{90}}+2 {q^{91}}+4 {q^{92}}+6 {q^{93}}
-2 {q^{95}}- 4 {q^{96}}-2 {q^{97}}
-2 {q^{98}}+{q^{99}}+4 {q^{100}}+{q^{101}}
\\ &
-{q^{103}}-2 {q^{105}}- {q^{109}}+2 {q^{111}}+2 {q^{113}}
-3 {q^{116}}-2 {q^{117}}+{q^{118}}+3 {q^{120}}+ {q^{121}}
\\ &
+3 {q^{122}}
-2 {q^{123}}-3 {q^{124}}-2 {q^{126}}+{q^{127}}+4 {q^{129}}- {q^{131}}
+{q^{132}}-{q^{133}}-3 {q^{135}}
\\ &
+3 {q^{136}}-{q^{137}}-{q^{138}}
+2 {q^{139}}+ {q^{140}}+{q^{141}}-3 {q^{142}}+3 {q^{143}}-2 {q^{144}}
-2 {q^{145}}
\\ &
+{q^{146}}+ {q^{147}}+{q^{148}}-3 {q^{149}}+4 {q^{150}}
-2 {q^{151}}-{q^{152}}+{q^{153}}+ {q^{154}}+{q^{155}}
-3 {q^{156}}
\\ &
+3 {q^{157}}-3 {q^{158}}-{q^{159}}+{q^{160}}+{q^{161}}+ {q^{162}}
-{q^{163}}+3 {q^{164}}-3 {q^{165}}-{q^{166}}-{q^{170}}
\\ &
+3 {q^{171}}- 2 {q^{172}}+{q^{174}}-{q^{177}}+{q^{178}}-2 {q^{179}}
+{q^{181}}+{q^{183}}+{q^{185}}- 2 {q^{186}}
\\ &
+{q^{188}}-{q^{189}}+{q^{190}}+{q^{192}}
-2 {q^{193}}+{q^{195}}-{q^{196}}+ {q^{197}}+{q^{199}}-2 {q^{200}}
+{q^{202}}
\\ &
-{q^{203}}+{q^{206}}-2 {q^{207}}+{q^{209}}+ {q^{213}}
-2 {q^{214}}+{q^{220}}-{q^{221}}-{q^{228}}+{q^{233}}-{q^{235}}
\\ &
+{q^{240}}- {q^{242}}+{q^{247}}-{q^{249}}+{q^{254}}-{q^{256}}
+{q^{261}}+{q^{268}}+{q^{275}}\big)
\end{align*}

\Addresses\recd

\end{document}